\documentclass[a4paper]{article}

\usepackage[francais]{babel}
\usepackage{delarray,verbatim,enumerate}
\usepackage{amsmath,amsthm,amstext,amsbsy,amssymb,amsfonts,amscd}
\usepackage[all]{xy}

\usepackage{ifpdf}
\ifpdf
	\usepackage{aeguill}
\else
	\usepackage[T1]{fontenc}
\fi

\newtheorem{Th}{Th\'eor\`eme}[]

\newtheorem{Prop}[Th]{Proposition}
\newtheorem{Cor}[Th]{Corollaire}

\newtheorem{Sco}[Th]{Scolie}
\newtheorem{Def} [Th]{D\'efinition}

\def\Preuve{\smallskip\noindent {\it Preuve.~}}

\def\Remarque{\smallskip\noindent {\it Remarque.~}}

\def\Remerciements{\smallskip\noindent {\it Remerciements.~}}

\font\teneufm=eufm10
\font\seveneufm=eufm7
\font\fiveeufm=eufm5
\newfam\gothfam
\textfont\gothfam=\teneufm \scriptfont\gothfam=\seveneufm
\scriptscriptfont\gothfam=\fiveeufm

\def\R{\mathbb{R}}		\def\Q{\mathbb{Q}}		\def\Z{\mathbb{Z}}	\def\N{\mathbb{N}}
\def\Zl{\mathbb{Z}_\ell} 	\def\Ql{\mathbb{Q}_\ell}	\def\Rl{\mathbb{R}_\ell}
\def\T{\mathbb{T}}	\def\D{\mathbb{D}}	\def\U{\mathbb{U}}	\def\C{\mathbb{C}}

		\def\ov{\overline}
\def\Im{\operatorname{Im}}	\def\Hom{\operatorname{Hom}}	
\def\Ker{\operatorname{Ker}}	

\begin{document}

\title{\bf\Large Compactification $\ell$-adique de $\mathbb R$}

\author{ Jean-François {\sc Jaulent} }
\date{}
\maketitle
\bigskip

{\small
\noindent{\bf Résumé.} Nous construisons un groupe topologique compact $\Rl$, qui contient naturellement comme sous-groupes denses à la fois le groupe additif réel $\R$ et le groupe $\ell$-adique $\Ql$ pour un premier donné $\ell$, et nous en étudions quelques unes des  propriétés. En appendice, nous montons que cela donne une description arithmétique du solénoide $\ell$-adique $S_\ell$ classiquement défini en termes de feuilletages.}

\

{\small
\noindent{\bf Abstract.} For a given prime number $\ell$, we construct a compact topological group $\Rl$ which contains both the real additive group $\R$ and the $\ell$-adic one $\Ql$  as dense subgroups; thus we study some of its properties. This construction gives an arithmetic description of the so-called $\ell$-adic solenoid $S_\ell$ classically defined in terms of foliations.}\bigskip\bigskip

%%%%%%%%%%%%%%%%%%%%%%%%%%%%%%%%%%%%%%%%%%
\bigskip\bigskip\bigskip

\noindent{\large \bf Introduction}

\bigskip
%%%%%%%%%%%%%%%%%%%%%%%%%%%%%%%%%%%%%%%%%%

Soit $\ell$ un nombre premier et $R_\ell = \{0, 1, \dots , \ell-1 \}$ l'ensemble des entiers naturels inférieurs à $\ell$. Il est bien connu que toute série de la forme\smallskip

\centerline{$\sum^{+\infty}_{n_0} a_n \, \ell^n$,}\smallskip

\noindent où les $a_n$ sont pris dans $R_\ell$ et $n_0$ dans $\Z$, converge dans le complété $\ell$-adique $\Ql$ du corps des rationnels~; et qu'inversement tout nombre $\ell$-adique $q$ s'écrit de façon unique comme somme d'une telle série.\smallskip

De façon semblable, toute série de la forme\smallskip

\centerline{$\sum_{-\infty}^{n_0-1} a_n \, \ell^n$,}\smallskip

\noindent converge dans le corps des réels $\R$~; et inversement tout réel positif $r$ s'écrit comme somme d'une telle série, laquelle est même unique dès lors que les $a_n$ ne sont pas ultimement égaux à $\,\ell-1$.\smallskip

On peut donc demander s'il est possible de construire un espace topologique convenable, contenant de façon naturelle $\Ql$ ainsi que $\R$, équipé de surcroît de lois prolongeant celles de $\Ql$ et de $\R$, dans lequel donner un sens univoque à l'expression
$$
\sum^{+\infty}_{-\infty} a_n \, \ell^n,
$$
indépendamment de la césure utilisée pour la calculer en l'écrivant formellement~: \medskip

\centerline{$\sum_{-\infty}^{n_0-1} a_n \, \ell^n + \sum^{+\infty}_{n_0} a_n \, \ell^n$.}\bigskip

On s'aperçoit aisément qu'il n'est pas raisonnable d'espérer prolonger continûment la multiplication dans un tel espace, puisque les suites de terme général $\ell^n$ et $\ell^{-n}$ convergeant vers 0 dans $\Ql$ et dans $\R$ respectivement, leur produit, qui est constant et égal à 1, devrait alors converger également vers 0, ce qui est évidemment absurde. Il faut donc renoncer à obtenir un anneau topologique.\medskip

L'objet du présent travail est ainsi de construire très simplement un groupe topologique abélien et compact que nous notons $\Rl$, qui contient naturellement à la fois $\R$ et $\Ql$ comme sous-groupes denses, et qui permet de donner un sens précis à chaque somme $q+r$ d'un nombre $\ell$-adique $q$ et d'un nombre réel $r$, donc finalement à toute somme infinie du type précédent.\smallskip

Nous verrons en appendice que le groupe obtenu $\Rl$ s'identifie canoniquement au solénoide $\ell$-adique $S_\ell$ défini classiquement en termes de feuilletages et qui intervient naturellement dans des questions liées à l'étude de certains systèmes dynamiques ou de probabilités.

%%%%%%%%%%%%%%%%%%%%%%%%%%%%%%%%%%%%%%%%%%
\bigskip\bigskip

\noindent{\large \bf 1. Rappels sur $\R$ et sur $\Ql$ ; construction du groupe $\Rl$}

\bigskip
%%%%%%%%%%%%%%%%%%%%%%%%%%%%%%%%%%%%%%%%%%

Rappelons d'abord brièvement comment sont construits usuellement les complétés du corps $\Q$ des rationnels relativement à ses diverses valeurs absolues.\medskip

On sait, par le théorème d'Ostrowski (cf. e.g. \cite{Bo} ou \cite{La}), que les topologies de $\Q$ qui sont définies par des valeurs absolues non triviales sont~:
\begin{enumerate}
\item[(i)] la topologie réelle, qui correspond à la valeur absolue usuelle définie par~:\smallskip

\centerline{$|\,x\,|_\infty =\max\{x,-x\}$~;}
\item[(ii)] les topologies $\ell$-adiques, données par les valeurs absolues définies par~:\smallskip

\centerline{$|\,x\,|_\ell = \ell^{-\nu_\ell(x)}$,}
\end{enumerate}
\noindent où $\nu_\ell(x)$ désigne l'exposant du premier $\ell$ dans la factorisation du rationnel $x \ne 0$.
\medskip

Dans tous les cas, le procédé de complétion, valable d'ailleurs pour tout espace métrique $E$ (cf. e.g. \cite{Sc}), permet de regarder $\Q$ comme un sous-espace topologique dense d'un espace complet, qui est évidemment $\Q_\infty=\R$ dans le premier cas et $\Ql$ dans le second. Dans le contexte qui nous intéresse, la méthode algébrique standard consiste à prendre l'anneau des suites de Cauchy de rationnels (pour la métrique considérée) et à le quotienter par l'idéal maximal formé des suites qui convergent vers 0. On obtient ainsi un corps commutatif sur lequel la valeur absolue se prolonge de façon naturelle, ce qui en fait un espace métrique complet (et même localement compact) qui contient $\Q$ comme sous-espace dense. \medskip

Faisons choix maintenant d'un nombre premier arbitaire $\ell$. Nous avons alors~:

\setcounter{Th}{-1}
\begin{Th} Soit $\ell \in \N$ un nombre premier et $R_\ell = \{0, 1, \dots , \ell-1 \}$. Alors~:
\begin{enumerate}
\item[(i)] Tout réel $\,r$ s'écrit~: $r=\pm \sum_{-\infty}^{n_0-1} a_n \, \ell^n$, avec les $a_n$ dans $R_\ell$ et $n_0$ dans $\Z$~; cette écriture est unique si les $a_n$ ne sont pas ultimement égaux à $\ell-1$.
\item[(ii)] Tout nombre $\ell$-adique $\,q$ s'écrit de façon unique~: $q=\sum_{m_0}^{+\infty} a_n \, \ell^n$, avec les $a_n$ dans $R_\ell$~; et on a : $|\,q\,|_\ell=\ell^{-m_0}$, si $a_{m_0}$ est son premier chiffre non nul.
\item[(iii)] Dans chaque cas, les sommes finies $z=\pm \sum_{m_0}^{n_0} a_n \, \ell^n$ décrivent le sous-groupe $\Z[1/\ell]$, lequel est dense dans le groupe $\,\Ql$ comme dans $\R$.
\end{enumerate}
\end{Th}

Nous allons plonger $\R$ et $\Ql$ dans groupe topologique abélien en les recollant suivant le sous-groupe dense $\Z[1/\ell]$. Comme indiqué dans l'introduction, le point essentiel de notre construction est que toute somme finie $\sum_{n_0+1}^{m_0} a_n \ell^n$ doit avoir même image qu'elle soit regardée dans $\Ql$ ou dans $\R$. Cette observation nous amène donc à poser~:

\begin{Def}
{\rm Nous appelons} $\ell$-adifié du groupe additif $\R$ {\rm le groupe abélien $\Rl$ quotient de la somme directe du groupe $\ell$-adique $\Ql$ et du groupe additif réel $\R$ par l'image diagonale du sous-groupe $\Z[1/\ell]$ des rationnels $\ell$-entiers~:}\smallskip

\centerline{$\Rl = ( \Ql \oplus \R )/ \Z[1/\ell]$.}\smallskip

Le groupe $\Rl$ est un groupe divisible qui contient canoniquement $\,\Ql$ et $\R$.
\end{Def}

L'application naturelle de $\R$ (respectivement de $\Ql$) dans $\Rl$ étant évidemment injective, nous identifierons un réel $r$ et un nombre $\ell$-adique $q$ avec leurs images respectives dans $\Rl$, de sorte que nous avons ainsi par construction~:\smallskip

\centerline{$\Ql \cap \R = \Z[1/\ell]$.}\smallskip

De façon semblable nous noterons $\,q+r$ la somme dans le groupe abélien $\Rl$ du nombre $\ell$-adique $q$ et du réel $r$, somme qui n'est autre que la classe dans $\Rl$ du couple $(q,r)$ de $\Ql \oplus \R$. Un tel nombre sera dit un {\it réel-$\ell$-adique}.\smallskip

\Remarque Les groupes additifs $\Ql$ et $\R$ étant tous deux divisibles, le quotient $\Rl$ l'est aussi. On prendra garde cependant au fait que $\Rl$ n'est pas {\em uniquement} divisible~; ainsi on a bien, par exemple~:\smallskip

\centerline {$ q+r = n \, (q/n +r/n)$,}\smallskip

\noindent pour tout entier naturel $n > 1$~; mais (avec les mêmes conventions d'écriture)~:\smallskip

\centerline{$1/n + 0 \,\ne\, 0 +1/n$,}\smallskip

\noindent dès que $n$ est étranger avec $\ell$, puisque $1/n$ n'est pas alors dans l'anneau $\Z[1/\ell]$. Plus généralement, pour  $k = 0, \cdots, n-1$, les $n$ éléments $x_k$ donnés par :\smallskip

\centerline{$x_k = (q+k)/n + (r-k)/n$}\smallskip

\noindent sont alors distincts dans $\Rl$ mais vérifient tous l'identité : $n x_k = q+r$. \medskip

Il suit de là que le groupe $\Rl$ contient des éléments de torsion. Nous verrons plus loin que ceux-ci constituent même un sous-groupe dense de $\Rl$.\medskip

L'écriture $x = r+q$ n'étant pas unique du fait de notre construction, il est intéressant de disposer néanmoins d'une écriture {\it standard} d'un réel-$\ell$-adique $x$~:

\begin {Th}[Lemme de représentation]\label{LR}
Tout réel-$\ell$-adique $\,x$ s'écrit de façon unique :
$$
x = [x] + \{x\}, \quad avec \ [x] \in \Zl \ et \ \{x\} \in [0,1[.
$$
Nous disons que l'entier $\ell$-adique $[x]$ est la partie entière et que le réel $\{x\}$ est la partie fractionnaire de $x$. En d'autres termes, nous avons~: \smallskip

\centerline{$\Rl = \Zl \oplus [0,1[$.}
\end{Th}

\Preuve Partons d'un réel-$\ell$-adique $x = q+r$~; écrivons $q = \sum^{+\infty}_{n_0} a_n \ell^n$ la représentation canonique de sa partie $\ell$-adique $q$~; puis notons $q' = \sum^{-1}_{n_0} a_n \ell^n \in \N[1/\ell]$ la partie non entière de $q$ et $z = \sum^{+\infty}_{0} a_n \ell^n$ sa partie entière. Nous avons ainsi~: $x = z+r'$  avec $r'=r+q' \in \R$ et $z \in \Zl$. Ecrivant alors $x' = [x'] + \{x'\}$ la décomposition du réel $x'$ comme somme de sa partie entière $[x] \in \Z$ et de sa partie fractionnaire $\{x'\} \in [0,1[$, nous obtenons comme attendu~:
$$
x = (z+[x']) + \{x'\} \in \Zl + [0,1[,
$$
ce qui établit l'existence de la décomposition annoncée.\smallskip

L'unicité résulte immédiatement de  l'identité~:
$$
\Zl \cap [0,1[ = \Zl \cap \Z[1/\ell] \cap [0,1[ \,= \Z \cap [0,1[\, =\{0\}.
$$

\Remarque Dans la décomposition directe précédente, le sous-groupe réel $\R$ et le sous-groupe $\ell$-adique $\Ql$ s'écrivent respectivement~:
$$
\R = \Z \oplus [0,1[ \qquad {\rm et} \qquad \Ql = \Zl \oplus [0,1[_\ell,
$$
où $[0,1[_\ell$ désigne l'ensemble des rationnels $\ell$-entiers de l'intervalle $[0,1[$.

\begin{Cor}
Tout réel-$\ell$-adique $x$ s'écrit de façon unique comme somme \smallskip

\centerline{$\sum\limits^{+\infty}_{0} a_n\, \ell^n \, +\, \sum\limits_{-\infty}^{-1} a_n\,  \ell^n$,}\smallskip

\noindent d'un entier $\ell$-adique $q=\sum^{+\infty}_{0} a_n\, \ell^n \in \Zl$ et d'un réel $r=\sum_{-\infty}^{-1} a_n\, \ell^n \in [0,1[$, où les $a_n$ sont pris dans $R_\ell$ et non ultimement égaux à $\,\ell-1$ pour $n \ll 0$.
\end{Cor}

\begin{Cor}\label{La}
Soit $a \in \Z$ un entier relatif fixé.
Tout réel-$\ell$-adique $x$ s'écrit alors de façon unique~:

\centerline{$x = [x]_a + \{x\}_a, \quad avec \ [x]_a \in \ell^a \Zl \ et \ \{x\}_a \in [0,\ell^a[$.}\smallskip

En d'autres termes, nous avons la décomposition directe~: \smallskip

\centerline{$\Rl = \ell^a \Zl \oplus [0,\ell^a[$.}
\end{Cor}

%%%%%%%%%%%%%%%%%%%%%%%%%%%%%%%%%%%%%%%%%%
\bigskip\medskip

\noindent{\large \bf 2. Propriétés topologiques du groupe $\Rl$}

\bigskip
%%%%%%%%%%%%%%%%%%%%%%%%%%%%%%%%%%%%%%%%%%

Nous allons maintenant définir une distance sur $\Rl$ en faisant choix d'une valeur absolue $|\ \ |$ à valeurs dans $\R_+$. Pour éviter toute confusion, nous notons~:\smallskip

{\it (i)} $\ |r|_\infty = \max \{ r,-r \}$, la valeur absolue usuelle sur $\R$, et\smallskip

{\it (ii)} $|q|_\ell =  \ \ell^{-\nu _\ell(q)}$, la valeur absolue normalisée sur $\Ql$.\smallskip

\noindent Cela étant, la {\em valeur absolue} sur $\Rl$ peut être définie de la façon suivante~:

\begin{Def}
Nous appelons valeur absolue d'un réel-$\ell$-adique $x = q+r$ la quantité~:
$$
|q+r| = \inf_{z \in \Z[1/\ell]} \max\{\ |q+z|_\ell\ ,\  |r-z|_\infty \},
$$
laquelle est encore donnée par la formule~:
$$
|x| = \min \{ \max \{\ |[x]|_\ell\ ,\ |\{x\}|_\infty\}\ ,\  \max \{\ |[x]+1|_\ell\ ,\ |1-\{x\}|_\infty\}\}.
$$
\end{Def}

\Preuve Il s'agit de vérifier que les deux formules sont équivalentes, autrement dit que la borne inférieure dans la première formule est atteinte pour la valeur de $z$ qui correspond à la décomposition {\em standard}  $x=[x]+\{x\}$ donnée par la proposition 2, ou pour la même valeur augmentée de 1, ce qui correspond à la décomposition  $x=([x]+1)+(\{x\}-1)$.

Or l'existence même de la décomposition {\em standard} $x=[x]+\{x\}$ entraîne~:\smallskip

\centerline{$|x| \le  \max \{\ |[x]|_\ell\ ,\ |\{x\}|_\infty\} \le 1$~;}\smallskip

\noindent de sorte que la borne inférieure doit être recherchée parmi les décompositions $x = (q+z)+(r-z)$ qui vérifient simultanément $|q+z|_\ell \le 1$ et $|r-z|_\infty \le 1$, i.e. $q+z \in \Zl$ et $r-z \in [-1, +1]$~; ce qui ne laisse guère que les deux possibilités annoncées.

\Remarque La valeur absolue $|\ \ |$ ainsi définie coïncide donc avec la valeur absolue normalisée $| \ \ |_\ell$ sur $\Zl$ et avec la valeur absolue usuelle $| \ \ |_\infty$ sur $]\!-1/2,+1/2\,[$.

\begin{Th}
L'application $(x,y) \mapsto |x-y|$ est une distance sur $\Rl$ qui est invariante par translation~; elle fait de $\Rl$ un groupe topologique contenant à la fois $\Ql$ et $\R$ comme sous-groupes denses.
\end{Th}

\Preuve Vérifions d'abord que nous avons bien là une distance~: L'équivalence 

\centerline {$|x| = 0 \Leftrightarrow  x=0$}

\noindent étant immédiate, compte tenu de l'expression de la valeur absolue, ainsi que l'identité $|x| =|-x|$, seule pose problème l'inégalité triangulaire. Or, si $x = q+r$ est une décomposition d'un réel-$\ell$-adique $x$ donnant sa valeur absolue (en ce sens qu'on a~: $|x| = \max\{\ |q|_\ell\ ,\  |r|_\infty \}$) et $x' = q'+r'$ une décomposition analogue pour $x'$, il vient~:\smallskip

\centerline{$|x+x'| \le \max\{\,|q+q'|_\ell, |r+r'|_\infty \} \le \max\{\,|q|_\ell, |r|_\infty \} + \max\{\,|q'|_\ell , |r'|_\infty \}$~;}\smallskip

\noindent ce qui donne, comme attendu~:\smallskip

\centerline{$|x+x'| \le |x| + |x'| $ .}\medskip

En résumé, il suit de là que $\Rl$ est bien un groupe topologique.\smallskip

Ce point acquis, l'invariance par translation provient directement de la défi\-nition même de la distance.
Enfin, le sous-groupe $\Z$ étant dense dans $\Zl$ et l'ensemble $[-1/2,1/2[_\ell\, := \Z[1/\ell] \cap [-1/2,1/2[$ l'étant dans $[-1/2,1/2[$, il résulte de la remarque précédant le théorème que le sous-groupe\smallskip

\centerline{$\Z[1/\ell] = \Z \oplus [-1/2,1/2[_\ell$}

\noindent est lui-même dense dans~:

 \centerline{$\Rl = \Z_\ell \oplus [-1/2,1/2[$.}\smallskip

\noindent En particulier, il suit comme annoncé que $\Ql$ et $\R$ sont denses dans $\Rl$.\medskip

Nous pouvons maintenant établir la compacité du groupe $\Rl$, qui justifie le titre de cette note~:

\begin{Th} Le groupe topologique $\Rl$ est compact et connexe. Nous disons que c'est le compactifié $\ell$-adique du groupe additif réel $\R$.
\end{Th}

\Preuve Pour voir que $\Rl$ est compact, considérons la suite exacte courte~:

\begin{displaymath}
\xymatrix{0 \ar[r] &  \Zl \ar[r]^\iota & \Rl \ar[r]^{\!\!\!\!\!\!\!\!\!\!\!\!\pi} & \T = \R / \Z \ar[r] &  0}
\end{displaymath}

\noindent induite par l'injection canonique $\iota$ de $\Zl$ dans $\Rl$. D'un coté, l'expression de la valeur absolue donnée plus haut montre que $\iota$ est isométrique~; de l'autre, la décomposition canonique $\Rl = \Zl \oplus [0,1[$ montre que son conoyau s'identifie (comme groupe topologique) au tore $\T = \R / \Z$. La compacité de $\Rl$ résulte donc de celle des deux groupes $\Zl$ et $\T$.\par

Ce point acquis, la connexité est immédiate~: en effet, le sous-groupe dense $\R$ étant bien enchaîné (pour sa métrique naturelle donc, {\em a fortiori}, pour la métrique de $\Rl$), l'espace $\Rl$, qui est simultanément compact et bien enchaîné, est connexe.\medskip

Comme tout groupe compact (et, plus généralement, tout groupe localement compact), le groupe $\Rl$ admet donc une mesure de Haar~:

\begin{Cor} Le groupe $\Rl$ admet une unique mesure positive $\mu$ de masse 1 qui est invariante par translation~; celle-ci est
caractérisée sur les compacts élémentaires $\ell^n\Zl \oplus [\alpha, \beta]$ (avec $0 \le \alpha \le \beta < 1$) par la formule~:\smallskip

\centerline{$\mu(\ell^n\Zl \oplus [\alpha, \beta]\,)\,=\, \ell^{-n}\times (\beta  - \alpha)$.}\smallskip

En particulier les sous-groupes $\Ql$ et $\R$ sont $\mu$-négligeables dans $\Rl$.
\end{Cor}

\begin{Prop}Pour tout $x_0$ de $\Rl$, la composante connexe par arcs de $x_0$ dans $\Rl$ (i.e. l'ensemble des extrémités des chemins dans $\Rl$ d'origine $x_0$) est la droite affine $x_0 +\R\,$ passant par $x_0$. En particulier, l'espace $\Rl$ n'est pas connexe par arcs.
\end{Prop}

\Preuve Translatons par $1/2$ la décomposition donnée par le théorème 2~; puis écrivons~:

\centerline{$\Rl = \Zl \oplus [-1/2, 1/2\,[$,}\smallskip

\noindent et notons ici $x=z+r$ la décomposition correspondante d'un  réel-$\ell$-adique $x$. Par définition de la distance sur $\Rl$, pour tout $x$ de $\Rl$, la boule ouverte de centre $x$ et de rayon $1/\ell$ est caractérisée par~:\smallskip

\centerline{$B(x,1/\ell) = x \;+\; \ell \Zl \; +\;  ]\!-1/\ell, 1/\ell\,[$.}\smallskip

Soit alors $X \mid t \mapsto x(t)$ un chemin dans $\Rl$, d'origine $x_0 =x(0)$ et d'extrémité $x_1 = x(1)$. La continuité uniforme de $X$ sur le segment $I=[0,1]$ nous fournit un $\eta > 0$ tel que pour toute $\eta$-chaîne $0=t_0 < t_1 < \dots < t_m = 1$ de $[0,1]$, les points $x_k=x(t_k)$ vérifient les inégalités~: $|x_k - x_{k-1}| < 1/\ell$, pour $k=1, \dots ,m$.\par

En particulier, sur chaque intervalle $I_k = I \;\cap\; ]x_k-\eta , x_k+\eta[$ la fonction $x(t)$ prend ses valeurs dans la boule $B(x_k,1/\ell)$, ce qui permet d'écrire localement~: \smallskip

\centerline{$x(t) = x_k + z_k(t) + r_k(t)$, avec $z_k(t) \in \ell\Zl$ et $r_k(t) \in ]\!-1/\ell, 1/\ell\,[$.}\smallskip

\noindent celà étant, l'application $t \mapsto z_k(t)$, qui est continue de $I_k$ dans $\Zl$, est évidemment constante, puisque $I_k$ est connexe tandis que $\Zl$ est totalement discontinu. Il suit de là que $x(t)-x_k$ est réel pour $t \in I_k$~; et finalement que $x(t)-x_0$ est à valeurs dans $\R$ pour tout $t \in I$.\smallskip

Réciproquement, pour tout $x_1=x_0+r \in x_0+\R$, l'application $X \,| \,t \mapsto x(t) =x_0 + t\,r$ est clairement un chemin d'origine $x_0$ et d'extrémité $x_1$.\medskip

\Remarque Convenons de dire qu'une partie $C$ de $\Rl$ est {\em convexe} lorsque pour tout couple $(a,b)$ d'éléments de $C$ on a simultanément~:\smallskip

\centerline{$b-a \in \R$, \quad ainsi que l'inclusion \quad $[a,b] = \{a+t(b-a)\,|\,t \in [0,1]\} \subset C$.}\smallskip

\noindent Avec ces conventions, la composante connexe par arcs d'un élément $x$ de $\Rl$ est tout simplement le plus grand convexe contenant $x$.

%%%%%%%%%%%%%%%%%%%%%%%%%%%%%%%%%%%%%%%%%%
\bigskip\bigskip

\noindent{\large \bf 3. Sous-groupe de torsion $\Rl^{\rm tor}$ et sous-groupes fermés de $\Rl$}

\bigskip
%%%%%%%%%%%%%%%%%%%%%%%%%%%%%%%%%%%%%%%%%%

Précisons d'abord la structure du sous-groupe de torsion de $\Rl$~:

\begin{Prop}
L'ensemble $\Rl^{\rm tor}$ des éléments de torsion de $\Rl$ est un sous-groupe dénombrable dense de $\Rl$ et s'identifie au quotient $\Q / \Z[1/\ell]$ via le morphisme $\tau$ qui à un rationnel $s\in \Q$ associe la classe $\tau(s)$ dans $\Rl$ du couple $(s,-s) \in \Ql \oplus \R$.
\end{Prop}

En particulier, il vient :

\begin{Cor} Pour tout entier naturel $m>0$ étranger à $\ell$, le groupe $\Rl$ possède un unique sous-groupe de cardinal $m$~: il est cyclique et c'est le noyau $_{m\!}\Rl$ de la multiplication par $m$~; et c'est aussi l'image par $\tau$ du sous-groupe $(1/m)\Z[1/\ell]$ de~$\Q$.
\end{Cor}

\Preuve Soit $x=q+r$ un élément de $m$-torsion dans $\Rl$, i.e. vérifiant $mx=0$, pour un $m>0$ de $\N$. De l'égalité $mx=mq+mr=0$ dans $\Rl$, nous concluons~:\smallskip

\centerline{$mq =-mr \in \Z[1/\ell]$,} 

\noindent disons :

\centerline{$mq=-mr= n/\ell^a$ avec $n \in \Z$ et $a\in \N$~;}\smallskip

\noindent ce qui montre que $q=-r$ est un nombre rationnel. En d'autres termes, l'application~:\smallskip

\centerline{$\tau \quad | \quad \Q \ni s \mapsto s-s \in \Rl$,}\smallskip

\noindent qui à un rationnel $s$ associe la classe du couple $(s,-s)$ dans $\Rl$, est un épimorphisme  du groupe additif $\Q$ sur le sous-groupe de torsion $\Rl^{\rm tor}$. Comme son noyau est clairement $\Z[1/\ell]$, nous obtenons par passage au quotient l'isomor\-phisme annoncé~:\smallskip

\centerline{$\Rl^{\rm tor} \simeq \Q / \Z[1/\ell]$.}\smallskip

\noindent En particulier, pour tout $m$ étranger à $\ell$, le sous-groupe de $m$-torsion de $\Rl$ est le groupe cyclique~: 

\centerline {$_{m\!}\Rl = \tau((1/m)\Z[1/\ell]) \simeq (1/m)\Z[1/\ell]/ \Z[1/\ell] \simeq \Z/m \Z$.}
\medskip

Enfin, la densité de $\Rl^{\rm tor}$ dans $\Rl$ résulte de celle de $\Q$ dans le produit $\Ql \times \R$.

\begin{Cor} Le sous-groupe de torsion $\Rl^{\rm tor}$ ne possédant pas d'élément d'ordre $\ell$, il suit en particulier que $\Rl$ lui-même est {\em uniquement} $\ell$-divisible, i.e. que l'on~a~:\smallskip

\centerline{$\forall x \in \Rl \quad \exists ! \; y \in \Rl \quad x=\ell y$.}\smallskip

\noindent On retrouve là le fait que $\Rl$ est naturellement un module sur l'anneau $\Z[1/\ell]$.
\end{Cor}

Intéressons nous maintenant aux sous-groupes {\em compacts} de $\Rl$.\smallskip

Il est bien connu que les sous-groupes fermés $G$ de $\R$ autres que $\R$ lui-même, sont discrets et de la forme $a\Z$ pour un $a$ de $G$. Il en résulte que les sous-groupes fermés du tore $\T = \R / \Z$, autres que $\T$ lui-même, sont les sous-groupes cycliques discrets $\T_m = (1/m)\Z/ \Z \simeq \Z/m\Z$, pour $m\ge 1$. \par

D'un autre côté, on sait que les sous-groupes fermés de $\Ql$, autres que $\Ql$ lui-même, sont le sous-groupe nul, qu'il est commode d'écrire $\ell^\infty\Zl$, et les groupes compacts $\ell^a \Zl$, pour $a\in\Z$, lesquels sont isomorphes à $\Zl$. Ces sous-groupes ne sont plus {\em algébriquement} monogènes, comme dans le cas réel (ils ne sont pas isomorphes à $\Z$), mais le restent {\em topologiquement} (car fermetures respectives de sous-groupes isomorphes à $\Z$). \smallskip

\begin{Th} Les sous-groupes compacts du groupe topologique $\Rl$ sont~:
\begin{itemize}
\item[(i)] Les sous-groupe triviaux~: le sous-groupe nul $\,0$ et le groupe $\Rl$ lui-même.
\item[(ii)] Les sous-groupes finis $_{m\!}\Rl$ de $\Rl^{\rm tor}$, cycliques d'ordre $m$ étranger à $\ell$.
\item[(iii)] Les sous-groupes compacts $\ell^a\Zl$ de $\Ql$, procycliques et isomorphes à $\Zl$.
\item[(iv)]Les produits $_{m\!}\Rl \times \ell^a\Zl$, avec $\ell \nmid m$ et $a \in \Z$, lesquels sont procycliques donc topologiquement monogènes.
\end{itemize}
\end{Th}

\Preuve Soit donc $G$ un sous-groupe fermé non trivial de $\Rl$ et $G_\ell = G \cap \Ql$ son intersection avec $\Ql$.\par
 Le sous-groupe $G_\ell$ n'est pas dense de $\Ql$, puisqu'il serait alors dense dans $\Rl$, ce qui entraînerait $G=\Rl$, contrairement à l'hypothèse. C'est donc que nous avons $G_\ell \subset \ell^b\Zl$ pour un $b$ de $\Z$. Et comme la topologie sur le groupe compact $\ell^b\Zl$ est induite aussi bien par celle de $\Ql$ que par celle de $\Rl$, il suit de là que $G_\ell$ est un sous-groupe compact de $\ell^b\Zl$. En fin de compte, nous avons donc soit $G_\ell=\{0\}$, soit $G_\ell=\ell^a\Zl$ pour un $a$ de $\Z$. Examinons successivement ces deux éventualités~:\smallskip

$(i)$ Si $G_\ell$ est nul, le groupe $G$ rencontre trivialement $\Zl$ de sorte que la restriction à $G$ 
de la projection canonique $\pi$ de $\Rl$ sur $\T = \Rl / \Zl$ est un isomorphisme de $G$ sur $\pi(G)$. Le groupe $G$  s'identifie alors à un sous-groupe fermé de $\T$, c'est à dire soit à $\T$ lui-même, soit à l'un de ses sous-groupes finis $\T_m$. Mais, si $\pi(G)$ était égal à $\T$, la représentation standard des éléments de $G$ ferait intervenir toutes les parties fractionnaires possibles~; de sorte que $G$ contiendrait en particulier un élément de la forme $z + 1/\ell$ avec $z \in \Zl$ et $1/\ell \in [0,1[$, que nous pouvons aussi bien écrire $(z + 1/\ell)+0$ avec $z +1/\ell \in \Ql \setminus \{0\}$, puisque $1/\ell$ est dans $\Z[1/\ell]$, ce qui est exclu, $G$ étant supposé renconter trivialement $\Ql$. Il suit donc~: $\pi(G) = \T_m$ pour un $m \ge 1$, comme annoncé. Ainsi $G$ est alors fini, et c'est le sous-groupe de $m$-torsion $_{m\!}\Rl$.\smallskip

$(ii)$ Si $G_\ell$ n'est pas nul, nous avons $G_\ell = \ell^a\Zl$ pour un $a$ de $\Z$~; écrivons~:\smallskip

\centerline{$\Rl = \ell^a \Zl \oplus [0,\ell^a[$,}\smallskip

\noindent et $\pi_a$ la projection canonique de $\Rl$ sur le quotient $\Rl/\ell^a\Zl \simeq \R / \ell^a\Z \simeq \T$. Ici encore, l'image $\pi_a(G)$ est un sous-groupe compact de $\T$, et ce n'est pas $\T$ car le groupe $G$ contiendrait alors en particulier un élément de la forme $z+\ell^{-(a+1)}$, pour un $z$ de $\ell^a\Zl$, c'est à dire $(z+\ell^{-(a+1)})+0$, avec $z+\ell^{-(a+1)} \in G_\ell \setminus \ell^a\Zl$, ce qui est absurde. Ainsi donc, l'image $\pi_a(G)$ est cyclique d'ordre, disons, $m$~; de sorte que dans la décomposition ci-dessus de $\Rl$, le sous-groupe $G$ est donné par~:

\centerline{$G = \ell^a \Zl \, \oplus \, \{0, \ell^a/m, \dots ,(m-1)\ell^a/m \}$.}\smallskip

En fin de compte, $G$ est alors le produit direct du sous-groupe compact $\ell^a\Zl$ de $\Ql$ et du sous-groupe fini $\, \ell^a({}_{m\!}\Rl) = {}_{m\!}\Rl$. Il suit de là que $m$ n'est pas un multiple de $\ell$ et que $G$ n'est autre que la préimage de $\ell^a\Zl$ par la multiplication par $m$~:

\centerline{$G = \frac{1}{m}\,\ell^a\Zl : = \{x\in\Rl \mid mx \in \ell^a\Zl\}$.} \medskip

En particulier~:

\begin{Cor}
Les sous-groupes fermés infinis de $\Rl$ sont $\Rl$ et les préimages de $\Zl$ dans les multiplications par les éléments de $\Z[1/\ell]$.
\end{Cor}

\begin{Cor} Les sous-groupes discrets de $\Rl$ sont les sous-groupes finis $_m\Rl$.
\end{Cor}

\Preuve Un tel sous-groupe est fermé donc compact et, par conséquent, fini.

\begin{Cor}
Soit $G = \Z \,x$ un sous-groupe monogène non nul de $\Rl$. Alors~:
\begin{itemize}
\item[(i)] Si $x$ est d'ordre $m$ dans $\Rl^{\rm tor}$, $G$ est le groupe cyclique discret ${}_{m\!}\Rl$.
\item[(ii)] Si $x$ s'écrit $\zeta+q$ avec $\zeta$ d'ordre $m$ dans $\,\Rl^{\rm tor}$ et $q \ne 0$ dans $\Ql$, de $\ell$-valuation $a \in \Z$, le goupe $G$ est dense dans le produit direct du groupe cyclique discret ${}_{m\!}\Rl$ et du groupe compact infini $\ell^a\Zl$.
\item[(iii)] Enfin, si $x$ n'est pas dans $\Rl^{\rm tor} + \Ql$, le sous-groupe $G$ est dense dans $\Rl$.
\end{itemize} 
\end{Cor}

\Preuve Cela résulte de la classification des sous-groupes fermés de $\Rl$ donnée par le Théorème~: si $G$ n'est pas dense dans $\Rl$, sa fermeture $\ov G$ est alors un sous-groupe propre compact de $\Rl$, ce qui montre que $x$ s'écrit alors $\zeta+q$, avec $\zeta$ dans $\Rl^{\rm tor}$ et $q$ dans $\Ql$. Soit alors $m$ l'ordre de $\zeta$ et $a \in \ov\Z$ la $\ell$-valuation de $q$.
\begin{itemize}
\item[(i)] Si $q$ est nul, $G$ est évidemment le groupe cyclique discret ${}_{m\!}\Rl$.
\item[(ii)] Sinon, $G$ contient $\Z\, mq$, qui est dense dans $\ell^a\Zl$, d'où le résultat annoncé.
\end{itemize}

%%%%%%%%%%%%%%%%%%%%%%%%%%%%%%%%%%%%%%%%%%
\bigskip\medskip

\noindent{\large \bf 4. Endomorphismes de $\Rl^{\rm tor}$ et dualité de Pontrjagin.}

\bigskip
%%%%%%%%%%%%%%%%%%%%%%%%%%%%%%%%%%%%%%%%%%

Considérons maintenant le dual de Pontrjagin du groupe $\Rl$, i.e. le groupe\smallskip

\centerline{$\widehat \R_\ell = \Hom_{\rm cont}(\Rl,\T)$}\smallskip

\noindent des morphismes continus de $\Rl$ dans le tore $\T=\R/ \Z \simeq\Rl / \Zl$.\smallskip

Si $\chi \in \widehat\R_\ell$ est un caractère continu sur $\Rl$, son image $\Im \chi$ est alors un sous-groupe compact du tore $\T$ et son noyau $\Ker \chi$ un sous-groupe compact de $\Rl$.\smallskip

Si $\chi$ n'était pas surjectif, son image $\Im \chi$ serait un sous-groupe fini $\T_m$ du tore $\T$, et il suivrait~:

\centerline{$(\Rl : \Ker \chi) = \mid \Im \chi \mid = m$,}\smallskip

\noindent contrairement à la classification des sous-groupes fermés de $\Rl$ donnée plus haut.\smallskip

Ainsi $\chi$ est surjectif, son image $\Im \chi = \T$ est connexe par arcs~; et comme $\Rl$ ne l'est pas, son noyau $\Ker \chi$ est non donc trivial. Or s'il était fini, nous aurions~:

\centerline{$\T = \Im \chi \simeq \Rl/{}_{m\!}\Rl \simeq \Rl$,}\smallskip

\noindent ce qui est encore exclu, pour des raisons de connexité par arcs. Il vient donc~:\smallskip

\centerline{$\Ker \chi =\frac{1}{m}\,\ell^a\Zl$, pour un $a$ de $\Z$ et un $m\ge 1$,}\smallskip

\noindent d'où, par factorisation, le diagramme commutatif où toutes les flêches sont des isomorphismes~:
\begin{displaymath}
\xymatrix{\Rl/ \frac{1}{m}\,\ell^a\Zl \ar[r]^{\bar\chi} \ar[d]_\phi & \Rl/\Zl \simeq \T \\
 \Rl/\Zl \simeq \T \ar[ru]&
}
\end{displaymath}
Dans celui-ci l'isomorphisme $\phi$ est induit par la multiplication par l'élément $m/\ell^a$ de $\Z[1/\ell]$~; et comme le groupe des automorphismes du tore $\T$ se réduit à $\{\pm 1\}$, il suit de là que le caractère $\chi$ est induit par la multiplication par $\pm m/\ell^a$.\smallskip

En d'autres termes, nous avons donc~:

\begin{Th}
Le groupe $\widehat \R_\ell = \Hom_{\rm cont}(\Rl,\T)$ des caractères continus de $\Rl$ s'identifie au groupe additif discret $\Z[1/\ell]$.
\end{Th}

\Remarque La suite exacte courte canonique de modules discrets~:
\begin{displaymath}
\xymatrix{0 \ar[r] &  \Z \ar[r]^{\!\!\!\!\!\!\!\hat\pi} & \Z[1/\ell] \ar[r]^{\!\!\!\!\!\!\!\!\!\!\!\!\hat\iota} & \T_{\ell^\infty} = \Z[1/\ell] / \Z \ar[r] &  0}
\end{displaymath}

\noindent peut ainsi être regardée comme duale de la suite exacte courte de modules compacts~:

\begin{displaymath}
\xymatrix{0 \ar[r] &  \Zl \ar[r]^\iota & \Rl \ar[r]^{\!\!\!\!\!\!\!\!\!\!\!\!\pi} & \T = \R / \Z \ar[r] &  0}.
\end{displaymath}

\begin{Cor}
Les endomorphismes continus non triviaux du groupe topologique $\Rl$ sont les applications $\phi$~:

\centerline{$x \mapsto \pm m/\ell^\nu \, x$\quad avec $\quad\pm m/\ell^\nu \in \Z[1/\ell]$.}\smallskip

\noindent En d'autres termes, ce sont les homothéties pour sa structure de $\Z[1/\ell]$-module.\par
En particulier, l'image de $\phi$ est $\,\Rl$ et son noyau le sous-groupe fini ${}_m\Rl$.
\end{Cor}

\Preuve Soit $\phi$ un endomorphisme continu non nul du groupe topologique $\Rl$. Puisque $\Rl$ est uniquement $\ell$-divisible, $\phi$ est en particulier un $\Z[1/\ell]$-morphisme.

Or, d'un côté, l'image $\phi(\R)$ du sous-groupe connexe par arcs $\R$ est alors un sous-groupe connexe par arcs (donc contenu dans $\R$) et il est non nul (sans quoi $\phi$ serait nul sur un sous-groupe dense)~; de sorte qu'il vient~: $\phi(\R) = \R$. Posant $\rho = \phi(1) \in \R$, nous concluons que $\phi(z)$ coïncide avec $\rho z$ pour tout $z \in \Z[1/\ell]$.

Et d'un autre côté, le compositum $\chi = \pi \circ  \phi$ de $\phi$ avec la surjection canonique $\pi$ de $\Rl$ sur $\T$ est alors un caractère continu de $\Rl$, donc de la forme~: $x \mapsto \rho_0 x$ modulo $\Zl$, pour un $\rho_0 \in \Z[1/\ell]$, en vertu du théorème ci-dessus. 

Rassemblant ces deux résultats, nous obtenons $(\rho-\rho_0)z \in \Zl$ pour tout $z \in \Z[1/\ell]$~; d'où, comme attendu~: $\rho = \rho_0 \in \Z[1/\ell]$. Et $\phi$, qui coïncide avec la multiplication par $\rho$ sur la partie dense $\Z[1/\ell]$, est l'homothétie de rapport $\rho$.\smallskip

\begin{Sco}
Avec les conventions d'écriture données dans la section 1, la préimage d'un réel-$\ell$-adique  $x=q+r$ de $\Rl$ (avec $q \in\Ql$ et $r \in \R$) par l'endomorphisme $\phi \; : \; z \mapsto \pm m/\ell^a \, z$ est formée des $m$ éléments :\smallskip

\centerline{$x_k = \pm \big((q\ell^{\nu}+k)/m + (r\ell^{\nu}-k)/m \big)$, pour $k=1,\cdots,m-1$.}

\end{Sco}

\Remarque  La construction du compactifié de $\R$ peut bien évidemment être menée en prenant simultanément plusieurs nombres premiers $\ell_1,\cdots,\ell_n$. 
L'exemple le plus éclairant eu égard à la numération usuelle,  est ainsi celui du compactifié décimal $\R_{10}$, {\it i.e.} du produit amalgamé de $\R$ et de $\Q_{10}=\Q_2\times\Q_5$ au-dessus de l'anneau des décimaux $\D=\Z[1/10]$. Dans ce cas, les endomorphismes continus du groupe $\R_{10}$ sont les homothéties pour sa structure de $\D$-module.

 Bien entendu, il est aussi possible de faire intervenir simultanément toutes les places de $\Q$. Dans ce cas, le groupe obtenu $\widehat\R$, qui est le produit amalgamé de $\R$ et de $\widehat\Q \simeq \prod_\ell \Q_\ell$ au-dessus de $\Q$\,, s'identifie au dual de Pontrjagin de $\Q$.

%%%%%%%%%%%%%%%%%%%%%%%%%%%%%%%%%%%%%%%%%%
\bigskip\medskip

\noindent{\large \bf 5. Appendice : lien avec le solénoïde $\ell$-adique}

\bigskip
%%%%%%%%%%%%%%%%%%%%%%%%%%%%%%%%%%%%%%%%%%

Pour étudier la dynamique d'une application continue $f$ sur un espace topologique $X$, il est classique d'introduire la limite du système projectif :
$$
Y = \varprojlim \; ( X \overset{f}\longleftarrow X \overset{f}\longleftarrow X \overset{f}\longleftarrow \cdots)
$$
défini par les itérées de $f$ ({\it cf. e.g.} \cite{CC}, Ch. 11). Dans le cas particulier où $X$ est le cercle unité $\U = \{z \in \C \mid |z| =1 \}$ du plan complexe, et $f$ l'application $z \mapsto z^\ell$, l'espace obtenu est le {\em solénoïde $\ell$-adique} : $$
S_\ell = \{(y_n)_{n \in \N} \in \U^\N \mid y_n =y_{n+1}^\ell \quad \forall n \in \N \}.
$$
Maintenant, en identifiant le tore $\T=\R/\Z$ au cercle unité $\U$ via l'application exponentielle $t \mapsto \exp(2i\pi t)$, on obtient l'isomorphisme de groupes topologiques :
$$
S_\ell \simeq \{(x_n)_{n \in \N} \in \T^\N \mid x_n = \ell x_{n+1} \quad \forall n \in \N \}.
$$
Et, sous cette dernière forme, il est aisé de reconnaitre le compactifié $\Rl$ de $\R$ défini plus haut :

\begin{Prop}
L'application qui à un réel-$\ell$-adique $x=\sum_{n=-\infty}^{+\infty}\; a_n \ell^n \in \Rl$ associe la famille des troncatures $(x_n)_{n\in\N}$, avec $x_n = \sum_{k=-\infty}^{n-1}\; a_k\ell^{n-k}+\Z \in \T$, composée avec l'exponentielle complexe $t \mapsto \exp(2i\pi t)$, réalise un isomorphisme de groupes topologiques $\varphi$ du compactifié $\Rl$ de $\R$ sur le solénoïde $\ell$-adique $S_\ell$.
$$
\varphi \; :\; \Rl \ni x=\sum_{n=-\infty}^{+\infty}\; a_n \ell^n \mapsto \Big( \exp \big( 2i\pi \sum_{k=-\infty}^{n-1}\; a_k \ell^{k-n}\big) \Big)_{n\in\N}\in S_\ell
$$
\end{Prop}

\Preuve Observons que $\varphi$ prend bien ses valeurs dans la limite projective $S_\ell \subset \U^\N$ et que c'est un morphisme de groupes. La condition de projectivité $x_n = \ell x_{n+1}$ montre que les chiffres $a_k$ (pour $k\in\Z$) associés à une famille cohérente $(x_n)_{n\in\N}$ de $\T^\N$ sont bien définis et qu'ils sont uniques ; en d'autres termes que $\varphi$ est bijective. Enfin, il est clair que $\varphi$ est continue pour les topologies naturelles de $\Rl$ et de $S_\ell$ , donc bicontinue puisque les espaces de départ et d'arrivée sont compacts, $\Rl$ d'après le théorème 7 et $S_\ell$ par le théorème de Tychonov.\smallskip

On retrouve ainsi le fait que le compactifié $\ell$-adique $\Rl$ est connexe mais non connexe par arcs, plus précisément qu'il est localement isomorphe comme tout solénoïde au produit d'un espace euclidien (en l'occurrence $\T$) et d'un espace totalement discontinu (en l'occurrence $\Zl$) ; ce qu'affirme explicitement le Lemme de représentation donné dans la section 1 (Th. 2).\medskip

\noindent {\it Coda}. Donnons pour conclure quelques illustrations simples de cette interprétation arithmétique du solénoïde $\ell$-adique :\smallskip

\begin{itemize}
\item[(i)] Transporté par $\varphi$ le plongement canonique $x\mapsto 0+x$ de $\R$ dans $\Rl$, introduit dans la section 1, devient le morphisme continu et injectif\smallskip

\centerline{$x \mapsto (\exp(  (2i\pi x/\ell^n))_{n\in\N}$,}\smallskip

\noindent de $\R$ dans $S_\ell$ usuellement considéré en analyse harmonique ({\em cf. e.g.} \cite{HR}) ou en théorie des probabilités ({\em cf. e.g.} \cite{Hy}).\smallskip

\item[(ii)]  L'étude algébrique des endomorphismes continus de $\Rl$ effectuée dans la section 4 montre directement que le  noyau $\Ker \phi$ d'un tel endomorphisme est toujours d'ordre étranger à $\ell$ ; ceci généralise notamment le résultat d'imparité donné dans \cite{AF} dans le cas $\ell=2$.\smallskip

\item[(iii)] La représentation explicite des {\em racines $n$-ièmes} d'un élément $\varphi(x)$ du solénoïde $p$-adique proposée dans \cite{BK} peut s'obtenir en transportant par $\varphi$ les points de $n$-division de $x$ dans $\Rl$ : 
\begin{itemize}
\item si le réel-$\ell$-adique $x$ s'écrit $x=q+r$ avec $q\in \Zl$ et $r \in \R$, et si $n$ est étranger à $\ell$, ce sont les images par $\varphi$ des $n$ points de $n$-division de $x$ définis par $x_k = (q+k)/n + (r-k)/n$, pour $k=1,\cdots,n-1$, avec les conventions d'écriture expliquées dans la section 1 ; 
\item si $n$ s'écrit $n=m\ell^\nu$ avec $\ell \nmid m$, ce sont les images des $m$ points distincts $x_k = (q\ell^{-\nu}+k)/m + (r\ell^{-\nu}-k)/m$, pour $k=1,\cdots,m-1$, puisque le groupe additif $\Rl$ est {\em uniquement} $\ell$-divisible.
\end{itemize}
\end{itemize}
\medskip

\Remerciements L'auteur remercie tout particulièrement C. Deninger (Muenster) et H.W. Lenstra (Leiden) pour leurs commentaires éclairants sur la version préliminaire de ce travail.

%%%%%%%%%%%%%%% BIBLIOGRAPHIE %%%%%%%%%%%%%%%%%

{\small
\def\refname 
 }

\bigskip\noindent
{\small
\begin{tabular}{l}
{Jean-François {\sc Jaulent}}\\
Institut de Mathématiques de Bordeaux \\
Université {\sc Bordeaux} 1 \\
351, cours de la libération\\
F-33405 Talence Cedex\\
email : jaulent@math.u-bordeaux1.fr 
\end{tabular}
}

 \end{document}